\documentclass[graybox]{svmult}
\usepackage{mathptmx}       
\usepackage{helvet}         
\usepackage{courier}        
\usepackage{type1cm}        
\usepackage{makeidx}         
\usepackage{graphicx}        
\usepackage{multicol}        
\usepackage[bottom]{footmisc}
\usepackage{amsmath, amsfonts, amssymb}
\usepackage{color}
\usepackage{tikz}

\usepackage{amssymb}
\usepackage{amsfonts}
\usepackage{amsmath}

\newcommand\beq{\begin{equation}}
\newcommand\eeq{\end{equation}}

\usepackage{dsfont}

\newcommand{\dt}{\Delta t}

\def\Qcal{\mathcal{Q}}
\def\Rcal{\mathcal{R}}

\def\Ucal{\mathcal{U}}

\renewcommand{\to}{\rightarrow}

\newcommand{\grad}{\boldsymbol{\nabla}}

\newcommand{\RR}{\mathbb{R}}\newcommand{\R}{\mathbb{R}}
\newcommand{\ZZ}{\mathbb{Z}}

\makeindex             
\begin{document}
\title*{A Finite-Volume 
discretization 
of 
%
viscoelastic Saint-Venant equations for FENE-P fluids
}
\titlerunning{Saint-Venant/FENE-P shallow flows}
\author{S\'ebastien Boyaval}
\institute{S\'ebastien Boyaval 
\at 
Laboratoire d'hydraulique Saint-Venant (Ecole des Ponts ParisTech -- EDF R\& D -- CEREMA) Universit\'e Paris-Est, EDF'lab 6 quai Watier 78401 Chatou Cedex France, 
\& INRIA Paris MATHERIALS 
\email{sebastien.boyaval@enpc.fr}}
\maketitle
\abstract{Saint-Venant equations can be generalized to account for a viscoelastic rheology in shallow flows.
A 
Finite-Volume discretization 
for the 1D Saint-Venant system generalized to 
Upper-Convected Maxwell (UCM) fluids
was proposed in [Bouchut \& Boyaval, 2013],
which preserved a physically-natural stability property (i.e. free-energy dissipation) of the full system.
It invoked a relaxation scheme of Suliciu type for the numerical computation of approximate solution to Riemann problems.
Here, the approach is extended to the 1D Saint-Venant system generalized to the finitely-extensible nonlinear elastic fluids of Peterlin (FENE-P).
We are currently not able to ensure all stability conditions a priori, but numerical simulations went smoothly in a practically useful range of parameters.
\keywords{Saint-Venant equations, FENE-P viscoelastic fluids, Finite-Volume, simple Riemann solver, Suliciu relaxation scheme
\\[5pt]
}
}
\section{Introduction}
\label{sec:intro}

Saint-Venant equations standardly model shallow free-surface gravity flows 
and can be generalized to account for the viscoelastic rheology of non-Newtonian fluids \cite{bouchut-boyaval-2015}, 
Upper-Convected Maxwell (UCM) fluids in particular \cite{bouchut-boyaval-2013}.
Here, we consider a generalized Saint-Venant (gSV) system for \emph{finitely-extensible nonlinear elastic} fluids with Peterlin closure (FENE-P fluids)
in Cartesian coordinates
\begin{eqnarray}
\label{eq:gsv1}
\partial_{t}h + \partial_{x}(h u) & = & 0
\\
\partial_{t}(h u) + \partial_{x}\left( h u^2 + gh^2/2 + h N \right) & = & 0
\label{eq:gsv2}
\\
\label{eq:gsv11}
\lambda\left( \partial_{t}\sigma_{xx} + u\partial_{x}\sigma_{xx} +2(\zeta-1)\sigma_{xx}\partial_{x}u \right) & = & 1-{\sigma_{xx}}/({1-(\sigma_{zz}+\sigma_{xx})/\ell})
\\
\label{eq:gsv12}
\lambda\left( \partial_{t}\sigma_{zz} + u\partial_{x}\sigma_{zz} +2(1-\zeta)\sigma_{zz}\partial_{x}u \right) & = & 1-{\sigma_{zz}}/({1-(\sigma_{zz}+\sigma_{xx})/\ell})
\end{eqnarray}
for 1D $\vec{e}_y$-translation invariant flow along 
$\vec{e}_x$ under a uniform gravity field $-g\vec{e}_z$ with 
\begin{itemize}
 \item mean flow depth $h(t,x)>0$ (in case of a non-rugous flat bottom),
 \item mean flow velocity $u(t,x)$ (for \emph{uniform} cross sections), and
 \item a normal-stress difference 
$ N = G(\sigma_{zz}-\sigma_{xx})/(1-(\sigma_{zz}+\sigma_{xx})/\ell) $
given by 
conformation variables $\sigma_{zz},\sigma_{xx}>0$ constrained by $0<\sigma_{zz}+\sigma_{xx}<\ell$, 
a relaxation time $\lambda\ge0$ and an elasticity modulus $G>0$. 
\end{itemize}
Note that 
(\ref{eq:gsv1}-\ref{eq:gsv2}-\ref{eq:gsv11}-\ref{eq:gsv12}) formally reduces to the standard viscous Saint-Venant system 
with viscosity $\nu\equiv2\lambda G\ge0$ when $\ell\to\infty$, $\lambda\to0$ and $G\lambda<\infty$. 
Moreover
we have used the quite general 
Gordon-Schowalter derivatives 
with slip parameter $\zeta\in[0,\frac12]$ 
constrained by the hyperbolicity of the system (\ref{eq:gsv1}-\ref{eq:gsv2}-\ref{eq:gsv11}-\ref{eq:gsv12}).
(This follows after an easy computation similar to \cite{2016arXiv161108491B}.)

In this work, we discuss a Finite-Volume method to solve (numerically) the Cauchy problem for the nonlinear hyperbolic 1D system (\ref{eq:gsv1}-\ref{eq:gsv2}-\ref{eq:gsv11}-\ref{eq:gsv12}).
Standardly, we need to consider \emph{weak} solutions (in fact, to (\ref{eq1}-\ref{eq2}-\ref{eq3}-\ref{eq4}), see below) plus \emph{admissibility} constraints
that are physically-meaningful dissipation rules formalizing the thermodynamics second principle close to an equilibrium \cite{dafermos-2000}.
Here, we consider 
the \emph{inequality} associated with the companion conservation law for the \emph{free-energy} 
$$
F = h \left( \frac{u^2}2 + \frac{gh}2 - \frac{G}{2(1-\zeta)}\left( \ell\log\left((\ell-(\sigma_{xx}+\sigma_{zz}))/(\ell-2)\right) + \log(\sigma_{xx}\sigma_{zz}) \right) \right)
$$
that is, 
on denoting the impulse by $P=gh^2/2+hN$,
\begin{multline}
\label{eq:entropy}
 -\frac{Gh}{2(1-\zeta)\lambda} \Bigg( 
  \sigma_{xx}^{-1} \left(1-\frac{\sigma_{xx}}{1-(\sigma_{zz}+\sigma_{xx})/\ell}\right)^2
+ \sigma_{zz}^{-1} \left(1-\frac{\sigma_{zz}}{1-(\sigma_{zz}+\sigma_{xx})/\ell}\right)^2 \Bigg)
\\
=: D
\ge
  \partial_{t}F + \partial_{x}\left(u(F + P)\right)
\end{multline}
where the left-hand-side is obviously non-positive on the admissibility domain
$$
\Ucal^\ell := \{0< h,0<\sigma_{xx},0<\sigma_{zz},\sigma_{xx}+\sigma_{zz}<\ell\} \,.
$$
Note that we do not consider the vacuum state $h=0$ as admissible here, see \cite{2016arXiv161108491B}.

\section{Finite-Volume discretization of FENE-P/Saint-Venant}
\label{sec:fv}

Piecewise-constant approximate solutions to the Cauchy problem on $(t,x)\in[0,T)\times\R$ for the gSV system can be defined by a Finite-Volume (FV) method.
With a view to preserving $\Ucal^\ell$ and the dissipation \eqref{eq:entropy} after discretization 
by a FV method,
we choose $q=(h,hu,h\sigma_{xx},h\sigma_{zz})$ as discretization variable.
Indeed, 
the free-energy functional $F$ is 
\emph{convex} on the convex domain $\Ucal^\ell\ni q$ 
(this follows after an easy computation from \cite[Lemma 1.3]{bouchut-2004})
while it is not convex in the variable $\left(h,hu,h\Pi,h\Sigma\right)$
whatever smooth invertible functions $\varpi,\varsigma$ are used for the 
reformulation of gSV 
\begin{eqnarray}
\label{eq1}
\partial_{t}h + \partial_{x}(h u) & = & 0
\\
\label{eq2}
\partial_{t}(h u) + \partial_{x}\left( h u^2 + \frac{gh^2}2 + h N \right) & = & 0
\\
\label{eq3}
\partial_{t}(h \Pi) + \partial_{x}(h u \Pi) & = & 
\frac{ h^{3-2\zeta} \varpi'(\sigma_{xx}h^{2(1-\zeta)})}\lambda \left(1-\frac{\sigma_{xx}}{1-\frac{\sigma_{zz}+\sigma_{xx}}\ell}\right)
\\
\label{eq4}
\partial_{t}(h \Sigma) + \partial_{x}(h u \Sigma) & = &
\frac{ h^{2\zeta-1} \varsigma'(\sigma_{zz} h^{2(\zeta-1)})}\lambda \left(1-\frac{\sigma_{zz}}{1-\frac{\sigma_{zz}+\sigma_{xx}}\ell}\right)
\end{eqnarray}
with $\Pi=\varpi(\sigma_{xx}h^{2(1-\zeta)})$, $\Sigma=\varsigma(\sigma_{zz}h^{2(\zeta-1)})$
(computations are similar to \cite[Appendix]{bouchut-boyaval-2013}).
In the sequel, we therefore discretize a 
quasilinear system
with 
source 
\begin{equation}
\label{eq:quasilinear}
 \partial_{t}q + A(q)\partial_{x}q = S(q) \,,
\end{equation}
%
which we recall is not ambiguous here (for those discontinuous solutions built using a Riemann solver, at least)
thanks to the dissipation rule \eqref{eq:entropy}, see 
\cite{lefloch-2002,berthon-coquel-lefloch-2011,2016arXiv161108491B}. 
%

\subsection{Splitting-in-time}

In cell 
$(x_{i-1/2},x_{i+1/2})$, $i\in\mathbb{Z}$, with volume $\Delta x_i=x_{i+1/2}-x_{i-1/2}>0$ and center $x_i={(x_{i-1/2}+x_{i+1/2})}/2$,
we approximate $q$ solution to \eqref{eq:quasilinear} on $\RR_{\ge0}\times\RR\ni(t,x)$ by
$$ q_i^{n+1}\approx\frac1{\Delta x_i}\int_{x_{i-1/2}}^{x_{i+1/2}}q(t,x)dx,\ i\in\mathbb{Z}, t\in(t^n,t^{n+1}] $$ 
on a time grid $0=t^0<t^1<\ldots<t^n<t^{n+1}<\ldots<t^{N}=T$ where $\dt^n=|t^{n+1}-t^n|$ will be chosen small enough 
compared with $\Delta x=\sup_{i\in\mathbb{Z}}\Delta x_i<\infty$ to ensure stability. 

More precisly,
having in mind the numerical approximation of
a (well-posed) Cauchy problem for \eqref{eq:quasilinear} on $\RR_{\ge0}\times\RR$ with initial condition $q(t\to0^+)=q^0\in L^\infty(\RR)$,
and therefore starting from approximations $ q_i^0\approx\frac1{\Delta x_i}\int_{x_{i-1/2}}^{x_{i+1/2}}q^0(x)dx$, $i\in\mathbb{Z}$,
we shall define
the cell values $q_i^n$ 
in two steps for each $n=1,\ldots,N$:\\
(i) an approximate solution 
to the \emph{homogeneous} gSV system (i.e. without the source term $S$) on $[t^n,t^{n+1})$
is first computed by an explicit 
three-point scheme
\begin{equation}
 \label{eq:fvscheme}
 q_i^{n+1/2}=q_i^n-\frac{\Delta t^n}{\Delta x_i}\left(F_l(q_i^n,q_{i+1}^n)-F_r(q_{i-1}^n,q_i^n)\right) \,,
\end{equation}
(ii) an approximate solution to the full gSV system 
on $(t^n,t^{n+1}]$ is next computed 
as 
\begin{equation}
 \label{eq:source}
 q_i^{n+1}=q_i^{n+1/2}+\Delta t^nS(q_i^{n+1}) \,.
\end{equation}
Then, the 
scheme is consistent 
with weak solutions 
of (\ref{eq:gsv1}--\ref{eq:gsv2}) equiv. (\ref{eq1}--\ref{eq2})
\begin{equation}
\label{fullscheme}
 q_i^{n+1}=q_i^n-\frac{\Delta t^n}{\Delta x_i}\left(F_l(q_i^n,q_{i+1}^n)-F_r(q_{i-1}^n,q_i^n)\right)+ \Delta t^n S(q_i^{n+1}) 
\end{equation}
provided the two first flux components for the conservative part $(h,hu)$ of the variable $q$ (actually solutions to conservation laws)
are conservative $F_{l,h}=F_{r,h}:=F_h$, $F_{l,hu}=F_{r,hu}:=F_{hu}$ and consistent $F_h(q,q)=hu|_q$, $F_{hu}(q,q)=(h u^2 + gh^2/2 + h N)|_q$ as usual, 
and with the conservative interpretation 
(\ref{eq3}--\ref{eq4}) of (\ref{eq:gsv11}--\ref{eq:gsv12})
insofar as we next define $F_l$ and $F_r$ using a \emph{simple} approximate Riemann solver \cite{harten-lax-vanleer-1983}  
for (\ref{eq1}--\ref{eq2}--\ref{eq3}--\ref{eq4}).

Moreover, 
with a view to preserving $\Ucal^\ell$ and a 
discrete version of \eqref{eq:entropy}
\begin{equation}
 \label{eq:freenergydissipation}
 F(q_i^{n+1/2})-F(q_i^{n})+\frac{\Delta t^n}{\Delta x_i}\left(G(q_i^n,q_{i+1}^n)-G(q_{i-1}^n,q_i^n))\right)\le 0
\end{equation}
for a 
numerical free-energy flux function consistent with $G(q,q)=u(F+P)|_q$ in \eqref{eq:entropy},
%
%
in the sequel, we shall discuss the relaxation technique introduced by Suliciu 
as simple Riemann solver in step (i), because it proved satisfying for other close systems \cite{bouchut-2003,bouchut-2004,bouchut-boyaval-2013}
equipped with an ``entropy'' 
convex in the discretization variable like $F$ here. 
%
In the end, for the full scheme \eqref{fullscheme}, a consistent free-energy dissipation
\begin{equation}
\label{eq:freenergydissipationdiscrete}
F(q_i^{n+1})-F(q_i^{n})+\frac{\Delta t^n}{\Delta x_i}\left(G(q_i^n,q_{i+1}^n)-G(q_{i-1}^n,q_i^n))\right)\le \Delta t^n D(q_i^{n+1}) 
\end{equation}
then holds 
insofar $h_i^{n+1/2}=h_i^{n+1}$, $u_i^{n+1/2}=u_i^{n+1}$ and the convexity of $F$ imply
\begin{equation}
 \label{eq:sourcedissipation}
 F(q_i^{n+1})-F(q_i^{n+1/2})\le \Delta t^n D(q_i^{n+1})\le 0 \,.
\end{equation}
\begin{proof}
On noting $h_i^{n+1/2}=h_i^{n+1}$, $u_i^{n+1/2}=u_i^{n+1}$ it suffices to show that
\begin{eqnarray*}
\lambda\left( \sigma_{xx,i}^{n+1} - \sigma_{xx,i}^{n} \right) / \dt^n & = & 1-{\sigma_{xx,i}^{n+1}}/({1-(\sigma_{zz,i}^{n+1}+\sigma_{xx,i}^{n+1})/\ell})
\\
\lambda\left( \sigma_{zz,i}^{n+1} - \sigma_{zz,i}^{n} \right) / \dt^n & = & 1-{\sigma_{zz,i}^{n+1}}/({1-(\sigma_{zz,i}^{n+1}+\sigma_{xx,i}^{n+1})/\ell})
\end{eqnarray*}
imply \eqref{eq:sourcedissipation}.
Now, this is obvious, on noting the convexity of $F|_{h,u}$ in $(\sigma_{xx},\sigma_{zz})$ and
$$
\grad_{(\sigma_{xx},\sigma_{zz})}F|_{h,u} \cdot S = D
$$
since $\grad_{(\sigma_{xx}h^{2(1-\zeta)},\sigma_{zz}h^{2(\zeta-1)})}F \cdot (h^{2(\zeta-1)} S_{h\sigma_{xx}}, h^{2(1-\zeta)} S_{h\sigma_{zz}}) = D $ by design.
\end{proof}

\subsection{Suliciu relaxation of the Riemann problem without source}

For all time ranges $t\in[t^n,t^{n+1})$, $n=0\ldots N-1$, 
let us now define at each interface $x_{i+\frac12}$, $i\in\ZZ$, between cells $i$ and $i+1$
the numerical flux functions $F_l$ and $F_r$ 
\begin{equation}
\label{eq:FlFr}
\begin{array}{l}
 F_l(q_l,q_r) = F_0(q_l)-\int_{-\infty}^0\Bigl(R(\xi,q_l,q_r)-q_l\Bigr)d\xi,\\
 F_r(q_l,q_r) = F_0(q_r)+\int_0^\infty\Bigl(R(\xi,q_l,q_r)-q_r\Bigr)d\xi.
\end{array}
\end{equation} 
invoking an approximate solution $R\left((x-x_{i+1/2})/(t-t^n),q_i^n,q_{i+1}^n\right)$ to the Riemann problem for \eqref{eq:quasilinear}
with initial condition $q_i^n 1_{x<0} + 1_{x>0} q_{i+1}^n$ at $t=t^n$, and {\it any} $F_0$.

In this work, we propose as approximate solution that given by Suliciu relaxation 
\begin{equation}
\label{relaxation}
R(\xi,q_l,q_r)=L \Rcal\left(\xi,\Qcal_l,\Qcal_r\right) 
\end{equation}
i.e. 
the projection (operator $L$) onto $q\equiv(h,hu,h\sigma_{xx},h\sigma_{zz})$ 
of the \emph{exact} solution $\Rcal\left(\xi,\Qcal_l,\Qcal_r\right)$ of the Riemann problem for the 
system with relaxed pressure
\begin{equation} 
\label{eq:gsvrelaxedvacuum} 
\left\lbrace
\begin{aligned}
\partial_{t}h + \partial_{x}(h u)
& = 0
\\
\partial_{t}(hu) + \partial_{x}(h u^2 + \pi)
& = 0
\\
\partial_{t} (\sigma_{xx}h^{2(1-\zeta)}) + u \partial_{x} (\sigma_{xx}h^{2(1-\zeta)})
& = 0
\\
\partial_{t} (\sigma_{zz}h^{2(\zeta-1)}) + u \partial_{x} (\sigma_{zz}h^{2(\zeta-1)})
& = 0
\\
\partial_{t}(h\pi) + \partial_{x}(hu\pi+uc^2)
& = 0
\\
\partial_{t}\left(h(u^2/2+\hat e)\right) + \partial_{x}\left(hu(u^2/2+\hat e)+u\pi\right)
& = 0
\\
\partial_{t}c + u \partial_{x}c
& = 0
\end{aligned}
\right.
\end{equation}
and initial condition given by ($o=l,r$)
\begin{equation}
\label{initial} 
\Qcal_o = \left(h_o,(hu)_o,h_o^{1-2\zeta}(h\sigma_{xx})_o,h_o^{2\zeta-3}(h\sigma_{zz})_o,h_oP(q_o),(hu)_o^2/2h_o+e(q_o),c_o\right) 
\end{equation}
where 
$c_o(q_l,q_r)$ are chosen so as to ensure stability, that is the dissipation rule \eqref{eq:freenergydissipation} here
(see below). 
Note that \eqref{eq:gsvrelaxedvacuum} is a hyperbolic system which fully decomposes into linearly degenerate eigenfields,
so $\Rcal$ has an analytic expression (see formulas in \cite{bouchut-2004,bouchut-boyaval-2013}).
Note also: the Riemann solver $R$ is consistent under the CFL condition
\begin{equation}
\label{eq:CFL}
\Delta t^n \le \frac12 \inf_{i\in\ZZ} \frac1{\Delta x_i} \max\left( u_i^n-c_l(q_i^n,q_{i+1}^n)/h_i^n , u_i^n+c_r(q_i^n,q_{i+1}^n)/h_{i+1}^n \right) \,.
\end{equation}

It remains to specify a choice of functions $c_l,c_r$ preserving 
$\Ucal^\ell$ and ensuring \eqref{eq:freenergydissipation}. 

Although it is not clear whether our construction allows one to approximate solutions on any time ranges $t\in[0,T)$,
since the series $\sum_n\Delta t^n$ may be bounded uniformly for all space-grid choice ($\sup_i|u_i^n|$ may grow unboundedly as $n\to\infty$),
specifying such $c_l,c_r$ fully defines a computable scheme.
In particular, \eqref{eq:freenergydissipationdiscrete} then implies that 
\eqref{eq:source} at step (ii) always has at least 
one solution $q_i^{n+1}\in\Ucal^\ell$ for any $\dt^n$ fixed at step (i). 
(This can be shown using Brouwer fixed-point theorem like in \cite{barrett-boyaval-2011}.)


Note however a difficulty 
here for FENE-P fluids with $c_l,c_r$.
Suliciu relaxation approach \eqref{eq:gsvrelaxedvacuum} was retained at step (i) because the solver often allows one 
to preserve invariant domains like $\Ucal^\ell$ \emph{and} a dissipation rule \eqref{eq:freenergydissipation} %
through well-chosen $c_l,c_r$, see e.g. \cite{bouchut-2003,bouchut-2004,bouchut-boyaval-2013}.
Indeed, on noting the exact Riemann solution to \eqref{eq:gsvrelaxedvacuum}, to get \eqref{eq:freenergydissipation} on choosing
$ G(q_l,q_r) = u \Bigl( h\bigl(\frac{u^2}{2}+\hat e\bigr) + \pi \Bigr)|_{\Rcal(0,q_l,q_r)} $, it is enough that
$\forall q_l,q_r\in\Ucal^\ell$
\begin{equation}
\label{eq:subcharacteristiccondition}
 q_\xi:=L\Rcal\left(\xi,\Qcal_l,\Qcal_r\right)\in\Ucal^\ell
 \text{ and }
 h_\xi^2 \partial_h|_{h^{2-2\zeta}\sigma_{xx},h^{2\zeta-2}\sigma_{zz}} P(q_\xi) \le c_\xi^2 \,,\ \forall \xi\in\RR
\end{equation}
using $c_\xi = c_l(q_l,q_r)$ if $\xi<= u^*$ and $c_\xi = c_r(q_l,q_r)$ if $\xi>u^*$ with $u^*:=\frac{c_lu_l+\pi_l+c_ru_r-\pi_r}{c_l+c_r}$.

One can easily propose $c_l,c_r$ satisfying the first condition 
in \eqref{eq:subcharacteristiccondition}, i.e. 
\begin{eqnarray}
\label{condl}
\frac{1}{h_l^*} & = & \frac{1}{h_l}\left(1+\frac{c_r(u_r-u_l)+\pi_l-\pi_r}{(c_l/h_l)(c_l+c_r)}\right) >0
\\
\label{condr}
\frac{1}{h_r^*} & = & \frac{1}{h_r}\left(1+\frac{c_l(u_r-u_l)+\pi_r-\pi_l}{(c_r/h_r)(c_l+c_r)}\right) >0
\end{eqnarray}
as usual for Saint-Venant systems, 
plus the 
admissibility conditions ($o=l/r$)
\begin{equation}
\label{cond1}
(h_o^*)^{2(1-\zeta)} (h_o)^{2(\zeta-1)} \sigma_{zz,o} + (h_o^*)^{2(\zeta-1)} (h_o)^{2(1-\zeta)} \sigma_{xx,o} < \ell
\end{equation}
for any $\sigma_{zz,o},\sigma_{xx,o}>0$ satisfying $\sigma_{zz,o}+\sigma_{xx,o} < \ell$
(FENE-P fluids, see below).
But the second condition 
is usually treated for $\phi_o:h \to h \sqrt{ \partial_h|_{h^{2-2\zeta}_o\sigma_{xx,o},h^{2\zeta-2}_o\sigma_{zz,o}} P }$
monotone. 
Unfortunately, a lengthy (but easy) computation shows that the latter 
is not monotone here,
so the standard method 
to choose $c_l,c_r$ a priori 
does not apply. 

\subsection{Choice of relaxation parameter}

%
Let us 
treat the first part of \eqref{eq:subcharacteristiccondition} 
as usual and define 
$c_o=\max(h_o\sqrt{\partial_hP(q_o)}:=h_oa_o,\tilde c_o)$, $o=l/r$
such that the functions $\tilde c_o(q_l,q_r)$ ensure (\ref{condl}--\ref{condr}) 
and \eqref{cond1}.

First, let us inspect (\ref{condl}--\ref{condr}) classically  following \cite[section3.3]{bouchut-klingenberg-waagan-2010}.
Denoting $a_lY_l=(u_l-u_r)_++\frac{(\pi_r-\pi_l)_+}{h_la_l+h_ra_r}\ge0$, $a_rY_r=(u_l-u_r)_++\frac{(\pi_l-\pi_r)_+}{h_la_l+h_ra_r}\ge0$ so
$\frac1{h_o^*}\ge\frac{1-{h_oa_oY_o}/{c_o}}{h_o}$,
it then holds $(h_o^*)^{-1}\ge(h_o)^{-1}y_o>0$
with $y_o:=1-\frac{Y_o}{1+\alpha_o Y_o}\in(\frac{\alpha_o-1}{\alpha_o},1]$
provided one chooses $\tilde c_o>0$ such that 
$c_o\ge h_oa_o(1+\alpha_o Y_o)$ for $\alpha_o>1$,
which yields $h_o^*\in(0,h_o/y_o]$ thus (\ref{condl}--\ref{condr}) in particular.

On the other hand, let us now inspect \eqref{cond1}, which
rewrites with $h_o^*>0$ 
\begin{equation}
\label{eq:range}
 w_o A_o + w_o^{-1} B_o < 1 
 \Leftrightarrow 
 2A_ow_o\in\left( {1-\sqrt{1-4A_oB_o}} , {1+\sqrt{1-4A_oB_o}} \right) \subset \RR_{>0}
\end{equation}
with $w_o=(h_o^*/h_o)^{2(1-\zeta)}$, $A_o=\sigma_{zz,o}/\ell$, $B_o=\sigma_{xx,o}/\ell$ positive such that $A_o+B_o<1$ (hence $A_oB_o\le A_o(1-A_o)\le\frac14$)
and $2(1-\zeta)\in[1,2]$. 
The upper-bound in \eqref{eq:range} is satisfied with 
$ \alpha_o = (w_o^+)^{\frac1{2(1-\zeta)}}/((w_o^+)^{\frac1{2(1-\zeta)}}-1) >1$, 
on noting 
\begin{equation}
\label{eq:upper}
(w_o^+)^{\frac1{2(1-\zeta)}}
:=\left(({1+\sqrt{1-4A_oB_o}})/({2A_o})\right)^{\frac1{2(1-\zeta)}}
\ge\text{\small $\frac{\alpha_o}{\alpha_o-1}$}
\ge1/{y_o}
\ge{h_o^*}/{h_o} \,.
\end{equation}
It remains to ensure the lower bound in \eqref{eq:range}.
Obviously, $w_o^-:=\frac{1-\sqrt{1-4A_oB_o}}{2A_o}<1$ so one only needs to inspect the case $h_o^*\le h_o$. 
Now, with $a_lW_l=(u_r-u_l)_++\frac{(\pi_l-\pi_r)_+}{h_la_l+h_ra_r}\ge0$, $a_rW_r=(u_r-u_l)_++\frac{(\pi_r-\pi_l)_+}{h_la_l+h_ra_r}\ge0$,
if $ c_o \ge h_o a_o W_o ({(w_o^-)^{-\frac1{2(1-\zeta)}}}-1)^{-1} $ then holds
$$
(w_o^-)^{\frac1{2(1-\zeta)}}\le\left(1+{a_oh_oW_o}/{c_o}\right)^{-1} \le {h_o^*}/{h_o} \,.
$$

In the end, we claim the following choices
\begin{eqnarray}
\label{cl}
c_l = h_l \max\left( a_l + \alpha_l \left((u_l-u_r)_++\frac{(\pi_r-\pi_l)_+}{h_la_l+h_ra_r}\right) ,
 \beta_l \left((u_r-u_l)_++\frac{(\pi_l-\pi_r)_+}{h_la_l+h_ra_r} \right)
 \right)
\\
c_r = h_r \max\left( a_r + \alpha_r \left((u_l-u_r)_++\frac{(\pi_l-\pi_r)_+}{h_la_l+h_ra_r}\right) ,
 \beta_r \left((u_r-u_l)_++\frac{(\pi_r-\pi_l)_+}{h_la_l+h_ra_r} \right)
 \right)
\label{cr}
\end{eqnarray}
satisfy 
simultaneously  
(\ref{condl}--\ref{condr}) 
and \eqref{cond1} 
in a compatible way 
with $a_o=\sqrt{\partial_hP(q_o)}$, 
$\alpha_o = \max(2,(w_o^+)^{\frac1{2(1-\zeta)}}/((w_o^+)^{\frac1{2(1-\zeta)}}-1))$, 
$\beta_o = (w_o^-)^{\frac1{2(1-\zeta)}}/(1-(w_o^-)^{\frac1{2(1-\zeta)}}))$, 
$w_o^-=\frac{\ell-\sqrt{\ell-4\sigma_{zz,o}\sigma_{xx,o}}}{2\sigma_{zz,o} }$, 
$w_o^+=\frac{\ell+\sqrt{\ell-4\sigma_{zz,o}\sigma_{xx,o}}}{2\sigma_{zz,o} }$, 
for $o=l/r$.
%
%
Moreover, note that we have chosen $\alpha_o$ 
such that all subcharacteristic conditions \eqref{eq:subcharacteristiccondition} are satisfied in the $\ell\to\infty$ limit,
hence also the free-energy dissipation \eqref{eq:freenergydissipationdiscrete}.
Indeed, $\phi_o$ 
is monotone in the 
$\ell\to\infty$ limit and one can then apply the standard method to choose $c_l,c_r$
\cite{bouchut-boyaval-2013}.



\section{Numerical illustrations}
\label{sec:num}

\begin{figure}[t]
\centering
\includegraphics[scale=.35]{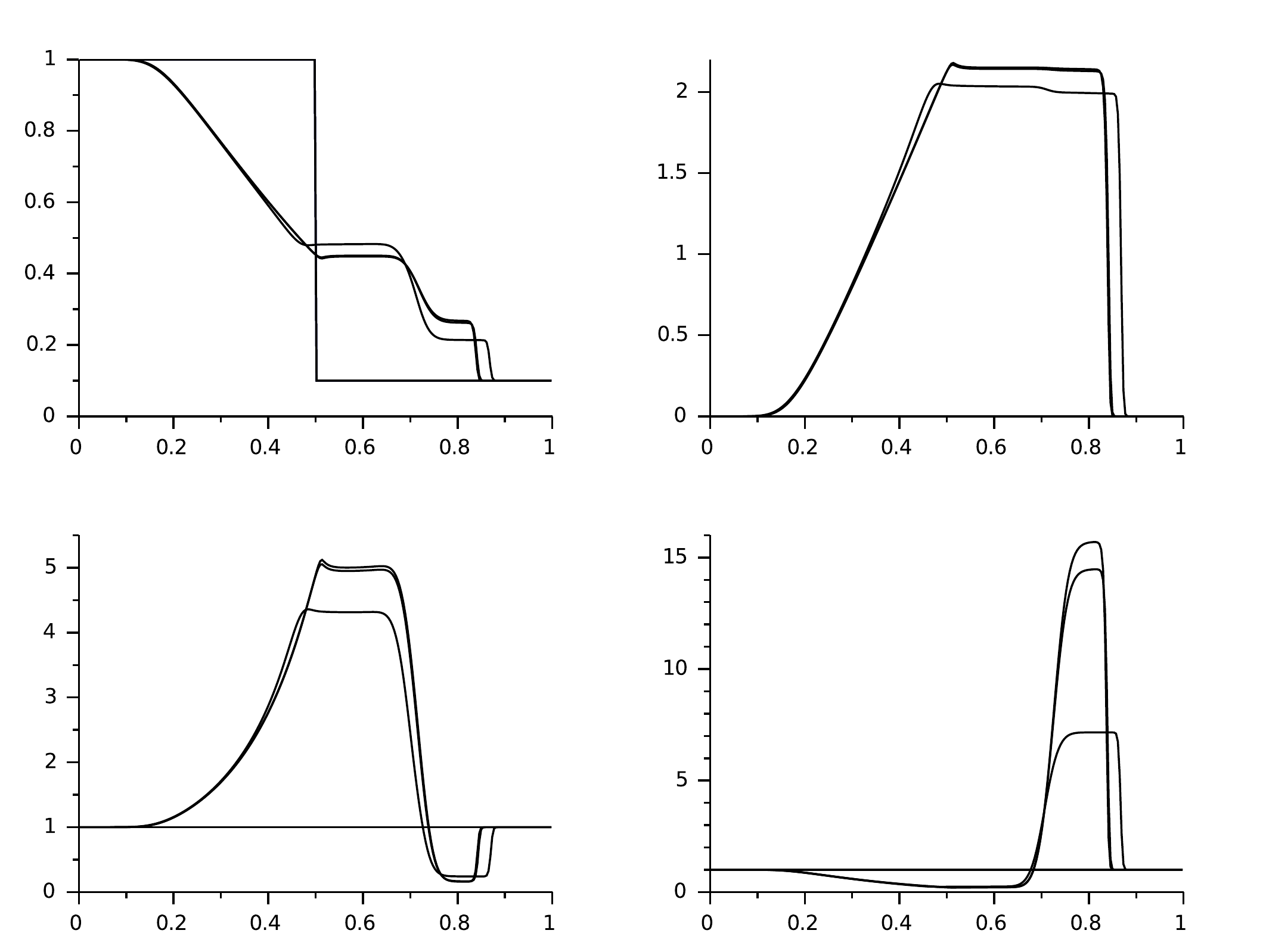}
\caption{Top: $h$ (left) and $u$ (right), bottom: $\sigma_{xx}$ and $\sigma_{zz}$. \label{fig:1}}
\end{figure}

We numerically approximate on $t\in[0,.1]$ the solution to a Riemann problem with
$$
\begin{cases}
 (h_l,u_l,\sigma_{xx,l},\sigma_{zz,l})=(1,0,1,1) & x<.5
 \\
 (h_r,u_r,\sigma_{xx,r},\sigma_{zz,r})=(.1,0,1,1) & x>.5 
\end{cases}
$$
as initial condition when $g=10$, $\zeta=0$, $G=.1$, 
$\lambda=.1$. 
In Fig.~\ref{fig:1}, we show the initial condition and the result at $t=.1$ when $\Delta x=2^{-8}$
for $\ell=10,100,1000$. 
Note the influence of the parameter $\ell$ on the stretch $\sigma_{xx}+\sigma_{zz}$.
On computing numerically the free-energy dissipation with the choice of relaxation parameter above, 
we have never observed the wrong sign, while the time-step did not 
go to zero.

\bibliographystyle{spmpsci}

\end{document}